\documentclass{article}

\usepackage{fancybox}
\usepackage{amsmath}
\usepackage{amssymb}
\usepackage[pdftex]{graphicx}
\usepackage{arydshln}

\newcommand{\vnorm}[1]{\left|\left|#1\right|\right|}

\newcommand{\IO}{\int_{\Omega}}

\newcommand{\eps}{\epsilon}
\newcommand{\bu}{{\boldsymbol{u}}}

\newcommand{\pu}{{\boldsymbol{p}}}
\newcommand{\zu}{{\boldsymbol{z}}}
\newcommand{\bnu}{{\boldsymbol{n}}}

\title{Wave-type threshold dynamics\\ and the hyperbolic mean curvature flow}
\author{Elliott Ginder\\ {\emph{Research Institute for Electronic Science}}\\{\emph{Hokkaido University}} \and Karel Svadlenka\\ {\emph{Department of Mathematics}}\\{\emph{Kyoto University}}}
\date{}
\begin{document}

\maketitle
\begin{abstract}
We introduce a method for computing interfacial motions governed by curvature dependent acceleration. Our method is a thresholding algorithm of the BMO-type which, instead of utilizing a diffusion process, thresholds evolution by the wave equation to obtain the desired interfacial dynamics. We also develop the numerical method and present results of its application, including an investigation of the volume preserving and multiphase motions.
\end{abstract}

\section{Introduction}

In this report, we propose a numerical scheme for the computation of the so-called hyperbolic mean curvature flow
\begin{equation}
\label{hmcf}
\frac{\partial^2 \gamma}{\partial t^2}(t,s) = -\kappa (t,s) \nu(t,s), 
\qquad \gamma(0,s)=\gamma_0(s), \quad \frac{\partial \gamma}{\partial t}(0,s)= v_0(s) \nu_0 .
\end{equation}
Here $\gamma : [0,T) \times I \to {\mathbb R}^2$ is a family of smooth curves, $\kappa$ denotes the curvature and $\nu$ is the unit outer normal vector.

This and similar types of hyperbolic problems were derived as models of oscillatory interface motions.
For example, \cite{[ls]} derives a model equation based on Hamilton's principle, 
considering stationary points of a geometrical action with local energy density consisting of kinetic and internal energy
$$ e = \frac{1}{2} \left( \Big| \frac{d}{dt} \gamma \Big|^2 + 1 \right) .$$
If the initial velocity is normal to the interface, the velocity remains normal during the motion, which leads to the equation
$$ \frac{d^2 \gamma}{dt^2} = e \kappa \nu - \nabla e . $$
Employing the variational structure, the authors show that the flow is locally well-posed and study conditions under which blow-up occurs.
In the case of graphs, they show uniqueness of solution satisfying certain entropy condition.

Similar equation 
$$ \frac{d^2 \gamma}{dt^2} = \kappa \nu - \nabla e , $$
which is found to be related to a hyperbolic Monge-Amp\`{e}re equation, 
is studied in \cite{[klw],[kw]} with an emphasis on the development of singularities.
The derivation of the equation is not given but its relation to equations for the motion of relativistic strings in Minkowski space is mentioned.

On the other hand, regarding damped oscillations, which appear on the solid-liquid interface of some crystals during melting or crystallization,
\cite{[gp],[rbn]} derived the model equation
$$ \rho \frac{dv}{dt}+ \beta v = (\psi + \psi'') \kappa - F, $$
(also see the references therein).
Here, $v, \frac{dv}{dt}$ are the normal velocity and normal acceleration, and $\psi, \beta, \rho, F$ denote the interfacial energy,
kinetic coefficient, effective density and a driving force for crystallization, respectively.

Equation (\ref{hmcf}) is also addressed in \cite{[hkl]}, where local existence and relations for the evolution of geometrical quantities are shown. 
It is a special case of a model equation for the motion of bubbles obtained in \cite{[kang]},
$$ \mu \frac{d u}{dt} = - p \nu - \sigma \kappa \nu + f - \mu u ( \nabla \cdot u- \nu \cdot Du \cdot \nu) -u \frac{d \mu}{dt} . $$
Here, $\frac{du}{dt}$ is the material derivative of the velocity vector, $\mu$ denotes mass density, $p$ is a factor corresponding to pressure, $\sigma$ is  surface tension and $f$ denotes additional sources of momentum.  

The numerical solution of the above hyperbolic mean curvature flows have been addressed only scarcely.
Except for straightforward front tracking schemes, \cite{[gp]} develops a crystalline algorithm for the motion of closed convex polygonal curves.
These methods cannot directly manage singularities, such as topological changes or the presence of junctions in the multiphase case, yet such problems can be resolved by adopting the level-set approach. A level-set method for curvature dependent accelerations based on the results of Sethian and Osher is presented in \cite{[kang]}. In this method, a nonlinear ill-posed problem has to be solved and it is not clear how the ideas can be extended to the multiphase setting.

Therefore, in this paper we aim at constructing a numerical scheme based on threshold dynamics of the BMO type.
The BMO algorithm was presented in \cite{[mbo]} and is based on the fact that frequent thresholding of the solution to the heat equation
approximates the evolution of level sets according to the standard (parabolic) mean curvature flow.
It is not obvious that such a method would be feasible for curvature accelerated motions. Nevertheless, investigating the formulae for the solution to the wave equation has revealed that level sets of the solutions whose initial condition is the signed distance to the interface, evolve with normal acceleration equal to their mean curvature, when the thresholding interval approaches zero.

In the following, we will present and formally justify the resulting algorithm (which we call HBMO, the hyperbolic BMO algorithm) together with its device for propagating velocities over the redistancing step without the need for explicit computation of velocities. We then comment on potential extensions of our scheme to more general flows, such as those involving phase volume constraints or the motion of junctions. We also include several numerical results and numerical confirmations regarding the proposed method.

\section{The HBMO algorithm}

The proposed HBMO algorithm for numerical approximation of the motion (\ref{hmcf}) in the case of a planar closed curve reads as follows. \\

\hrule

\vspace{0.2cm}

Given: initial curve $\gamma_0$, its normal velocity $v_0$, a final time $T$ and a time step $\tau =T/N$.
\begin{enumerate}
\item Extend $v_0$ in a suitable way to a neighborhood of $\gamma_0$.
\item For $t \in [0,\tau]$ solve the initial value problem
\begin{equation}
\label{wave_eq1}
u_{tt}(t,x) = \Delta u(t,x), \qquad u(0,x)=d_0(x), \qquad u_t(0,x)=-v_0(x), 
\end{equation}
where $d_0(x)$ is the signed distance function to $\gamma_0$.
\item Define $\gamma_1$ as the zero level set of $u(x,\tau)$.
\item For $n=1,2, \dots ,N-1$ repeat
\begin{enumerate}
\item For $t \in [0,\tau]$ solve the initial value problem
\begin{equation}
\label{wave_eq2}
u_{tt}(t,x) = 2\Delta u(t,x), \qquad u(0,x)=2d_{n}(x)-d_{n-1}(x), \qquad u_t(0,x)=0, 
\end{equation}
where $d_n(x)$ is the signed distance function to $\gamma_n$.
\item Define $\gamma_{n+1}$ as the zero level set of $u(x,\tau)$.
\end{enumerate}
\end{enumerate}

\hrule

\vspace{0.3cm}

In this section, we explain the ideas behind the derivation of the above HBMO algorithm and
provide a formal justification for its convergence. We divide the explanation into two parts. We begin by addressing the construction of the first curve $\gamma_1$ (steps 1 - 3 in the algorithm), which will clarify the reason why a threshold-type scheme for the wave equations leads to accelerations proportional to curvature. Then we treat further HBMO steps (step 4 in the algorithm), focusing on obtaining the propagation of interfacial velocities throughout the evolution.

\subsection{The first step of the HBMO algorithm}
\label{sec1BMO}

Given an initial closed smooth planar curve $\gamma_0$ and its smooth initial normal velocity $v_0$,
we construct an approximation $\gamma_1$ of the curve's evolution by HMCF at time $\tau$.
This is done as follows:
\begin{enumerate}
\item Extend $v_0$ to the neighborhood of $\gamma_0$ in ${\mathbb R}^2$ as
$$ v_0(x,y) = v_0(\hat{x},\hat{y}), $$
where $(\hat{x},\hat{y})$ is the orthogonal projection of $(x,y)$ on $\gamma_0$, i.e., the nearest point to $(x,y)$ on the curve $\gamma_0$.
\item Solve the initial value problem
\begin{equation}
\label{wave_eq}
u_{tt} = c^2 \Delta u, \qquad u(0,x)=d(x), \qquad u_t(0,x)=-v_0(x), 
\end{equation}
where $d(x)$ is the signed distance function to $\gamma_0$, and where
we take $c=1$ here, (but we keep the general coefficient $c$ for later convenience).  We remark that the velocity $v_0$ is defined only in a neighborhood of $\gamma_0$ but, since $\gamma_0$ is smooth and 
$\tau$ can be taken small enough, this fact will not hinder the subsequent step due to the finite speed of propagation.
\item Define $\gamma_1$ as the zero level set of $u(x,\tau)$.
\end{enumerate}

We now examine the evolution of the interface, making use of the explicit representation formula. The solution to (\ref{wave_eq}) is equal to $u(t,x)=v(ct,x)$, where $v$ is the solution of
$$ v_{tt} = \Delta v, \qquad v(0,x)=d(x), \qquad v_t(0,x)=-v_0(x)/c . $$
Here and in the sequel, $x=(x_1,x_2)$ and $y=(y_1,y_2)$ are points in $\mathbb{R}^2$.
Using the Poisson formula for $v$ and a change of variables,
we find that the formula for the original solution $u$ reads
\begin{equation}\label{transPoisson}
	u(t,x) = \frac{1}{2 \pi ct} \int_{B(x,ct)} \frac{d(y)+ \nabla d(y) \cdot (y-x) -tv_0(y)}{\sqrt{c^2t^2 - |y-x|^2}} \, dy.
\end{equation}

In order to analyze the resulting motion, 
we take a point on the interface and rotate and translate the coordinate system, so that the point becomes
the origin and the outer normal to the interface points in the direction of the $x_2$-axis.
We assume that the interface is smooth in a neighborhood of the origin. Then the value of the signed distance function at any point $y \in B(x,ct)$ can be approximated by the following Taylor expansion 
(see \cite{[essedoglu]}), provided that $x$ in (\ref{transPoisson}) is close to the origin and $t$ is sufficiently small
(depending on the smoothness of the interface):
$$ d(y_1,y_2) = y_2 + \frac{1}{2} \kappa y_1^2 +  \frac{1}{6} \kappa_{x_1}y_1^3 - \frac{1}{2} \kappa^2 y_1^2 y_2 + 
\sum_{| \alpha | = 4} e_{\alpha}(y) y^{\alpha} . $$
Here $\kappa$ is the curvature of the interface at the chosen point (now the origin) and the 
$e_{\alpha}$'s are smooth functions. We remark that the error functions, denoted by $e_{\alpha}$ (with $\alpha$ multiindex),
will vary from place to place and are always assumed to be smooth and bounded functions of their variables. 

We now proceed to calculate the contribution of each term in the signed distance expansion to the solution of wave equation $u(t,x)$.


For the first term $y_2$ we have
\begin{eqnarray*}
u^1(t,x) &=& \frac{1}{2 \pi ct} \int_{B(x,ct)} \frac{y_2 + \begin{pmatrix} 0 \\ 1 \end{pmatrix} \cdot \begin{pmatrix} y_1-x_1 \\ y_2- x_2 \end{pmatrix}}{\sqrt{c^2t^2 - |y-x|^2}} \, dy \\
&=& \frac{1}{2 \pi ct} \int_{B(0,1)} \frac{2(ct z_2+ x_2) - x_2}{ct \sqrt{1-|z|^2}} c^2t^2 \, dz \\
&=& \frac{1}{2 \pi} \int_{B(0,1)} \frac{x_2}{\sqrt{1-|z|^2}} \, dz \\
&=& \frac{x_2}{2 \pi} \int_0^1 \int_0^{2 \pi} \frac{r}{\sqrt{1-r^2}} \, d \theta \, dr \\
&=& x_2 
\end{eqnarray*}
Here we use the change of variables $z=(y-x)/ {ct}$ to transform the domain of integration to a fixed ball $B(0,1)$.
Also we use the fact that the function $z_2/\sqrt{1-|z|^2}$ is odd with respect to the $x_2$-axis, hence its integral over 
$B(0,1)$ vanishes.

By a similar calculation, the second term $\frac{1}{2} \kappa y_1^2$ gives
\begin{eqnarray*}
u^2(t,x) &=& \frac{1}{2 \pi c t} \frac{\kappa}{2} \int_{B(x,ct)} \frac{y_1^2 + \begin{pmatrix} 2y_1 \\ 0 \end{pmatrix} \cdot \begin{pmatrix} y_1 -x_1 \\ y_2- x_2 \end{pmatrix}}{\sqrt{c^2t^2 - |y-x|^2}} \, dy \\
&=& \frac{\kappa}{4 \pi} \int_{B(0,1)} \frac{3c^2t^2z_1^2+x_1^2}{\sqrt{1-|z|^2}} \, dz \\
&=& \frac{\kappa}{2} (c^2t^2+x_1^2) .
\end{eqnarray*}

The third term $\frac{1}{6} \kappa_{x_1}y_1^3$ contributes to the solution as follows:
\begin{eqnarray*}
u^3(t,x) &=& \frac{\kappa_{x_1}}{12 \pi c t} \int_{B(x,ct)} \frac{y_1^3 + \begin{pmatrix} 3y_1^2 \\ 0 \end{pmatrix} \cdot \begin{pmatrix} y_1 -x_1 \\ y_2- x_2 \end{pmatrix}}{\sqrt{c^2t^2 - |y-x|^2}} \, dy \\
&=& \frac{\kappa_{x_1}}{12 \pi} \int_{B(0,1)} \frac{9c^2t^2x_1z_1^2+x_1^3}{\sqrt{1-|z|^2}} \, dz \\
&=& \frac{\kappa_{x_1} x_1}{6} (3c^2t^2 + x_1^2).
\end{eqnarray*}

The fourth term $- \frac{1}{2} \kappa^2 y_1^2 y_2$ yields
\begin{eqnarray*}
u^4(t,x) &=& - \frac{\kappa^2}{4 \pi c t} \int_{B(x,ct)} \frac{y_1^2 y_2 + \begin{pmatrix} 2y_1y_2 \\ y_1^2 \end{pmatrix} \cdot \begin{pmatrix} y_1 -x_1 \\ y_2- x_2 \end{pmatrix}}{\sqrt{c^2t^2 - |y-x|^2}} \, dy \\
&=& -\frac{\kappa^2}{4 \pi} \int_{B(0,1)} \frac{3c^2t^2x_2z_1^2+6c^2t^2x_1z_1z_2 +x_1^2x_2}{\sqrt{1-|z|^2}} \, dz \\
&=& -\frac{\kappa^2x_2}{2} (c^2t^2 + x_1^2).
\end{eqnarray*}

The velocity term is approximated by its Taylor series as follows:
\begin{eqnarray*}
u^v(t,x) &=& \frac{1}{2 \pi ct} \int_{B(x,ct)} \frac{tv_0(y)}{\sqrt{c^2t^2 - |y-x|^2}} \, dy \\
&=& \frac{t}{2 \pi} \int_{B(0,1)} \frac{v_0(x+ctz)}{ \sqrt{1-|z|^2}} \, dz \\
&=& \frac{t}{2 \pi} \int_{B(0,1)} \frac{v_0(0)+\nabla v_0(0) \cdot (x+ctz) + 
\sum_{|\alpha| \geq 2} e_{\alpha}(t,x,z)(t,x)^{\alpha}}
{\sqrt{1-|z|^2}} \, dz \\
&=& v_0(0)t + \tfrac{\partial}{\partial x_1} v_0(0) t x_1 + 
\sum_{|\alpha| \geq 3} e_{\alpha}(t,x)(t,x)^{\alpha} ,
\end{eqnarray*}
since $\tfrac{\partial v_0}{\partial x_2}(0)=0$.

Finally, the error term in the signed distance expansion can be evaluated as
\begin{eqnarray*}
u^e(t,x) &=& \frac{1}{2 \pi ct} \sum_{| \alpha | = 4} \int_{B(x,ct)} \frac{ e_{\alpha}(y) y^{\alpha} + \nabla \left( e_{\alpha}(y) y^{\alpha}  \right) \cdot (y-x) }{\sqrt{c^2t^2 - |y-x|^2}} \, dy \\
&=& \frac{1}{2 \pi} \sum_{| \alpha | = 4} \int_{B(0,1)} \frac{e_{\alpha}(x+ctz) (x+ctz)^{\alpha} + \nabla \left( e_{\alpha}(y) y^{\alpha}  \right)|_{y=x+ctz} \cdot ctz}{ \sqrt{1-|z|^2}} \, dz \\
&=& \sum_{|\alpha| =4} e_{\alpha}(t,x)(t,x)^{\alpha} .
\end{eqnarray*}

The solution to the wave equation with initial condition $d$ in the neighborhood of the origin can thus be written in the following way:
\begin{eqnarray}
u(t,x) &=&  x_2 + \frac{\kappa}{2} (c^2t^2+x_1^2) + \frac{\kappa_{x_1} x_1}{6} (3c^2t^2 + x_1^2) -\frac{\kappa^2x_2}{2} (c^2t^2 + x_1^2) - u^v(x,t)+u^e(x,t) \nonumber \\
&=&  x_2 + \frac{\kappa}{2} (c^2t^2+x_1^2) + \frac{\kappa_{x_1} x_1}{6} (3c^2t^2 + x_1^2) -\frac{\kappa^2x_2}{2} (c^2t^2 + x_1^2) \nonumber \\
\label{sol1step}
&& - v_0(0)t - \tfrac{\partial v_0}{\partial x_1} (0) tx_1 + 
\sum_{|\alpha| \geq 3} e_{\alpha}(t,x)(t,x)^{\alpha}.
\end{eqnarray}

From this result we can make the following observation.
If the distance travelled by the interface in the normal direction (i.e., the $x_2$-direction) after time $\tau$ is denoted by $\delta_0$,
then this distance can be calculated from the relation $u(\tau,0,\delta_0)=0$. In particular,
$$  0=u(\tau,0, \delta_0) =  \delta_0 + \frac{1}{2} \kappa c^2 \tau^2 -\frac{1}{2}  \delta_0 \kappa^2 c^2 \tau^2 - v_0(0)\tau 
+ \sum_{|\alpha| \geq 3} e_{\alpha}(\tau,\delta_0)(\tau,\delta_0)^{\alpha}
. $$
This relation implies that, in terms of the order in $\tau$, the second order approximation of $\delta_0$ is $v_0(0)\tau-\frac{1}{2} \kappa c^2 \tau^2$,
since the error term is of order $O(\tau^3+\delta_0^3)$.
Therefore, the error term can be estimated by $O(\tau^3)$ and solving the above equation for $\delta_0$, we obtain
\begin{eqnarray}
\delta_0 =  \frac{v_0(0)\tau-\frac{1}{2} \kappa c^2\tau^2 +O(\tau^3)}{1- \frac{1}{2} \kappa^2 c^2\tau^2} 
= v_0(0)\tau - \frac{1}{2} \kappa c^2 \tau^2 +O(\tau^3) . 
\end{eqnarray}
This means that the interface moves with initial velocity $v_0$ and with an acceleration equal to $- c^2 \kappa$, the $c^2$-multiple of the curvature.
Thus taking $c=1$, we obtain the desired approximation of the interface $\gamma_1 \approx \gamma(\tau)$.

\subsection{Further steps of HBMO}

If we want to apply the above scheme to all time steps, the velocities along the interface would need to be computed. This would severely complicate the numerical solution, even if the velocity field was known. The difficulties inherent to such an approach are made clear, for example, when considering evolutions that involve topological changes (e.g., interfaces that split apart, or contact each other). On the other hand, if we reset the initial velocity to zero, the interface will always accelerate from zero velocity and the accumulated velocities will not propagate. Hence, we need to provide a modification of the above scheme which inherits the velocity from the previous time step and accelerates the interface depending on its curvature.

In order to design a scheme that avoids computation of interface velocities, 
it is necessary to account for the position of the interface at both the previous and the current time step.
Thus, it is natural to consider a vanishing initial velocity in the wave equation and an initial condition of the form
\begin{equation}
\label{initprev}
u(0,x) = a d_{n-1}(x) + bd_n(x).
\end{equation}
Here $d_n$ denotes the signed distance function to the present interface $\gamma_n$, $d_{n-1}$ is the signed distance function to the interface
at previous time step $\gamma_{n-1}$, and both $a$ and $b$ denote real numbers (see figure \ref{xi-coords} below).
\begin{figure}[ht!]
\begin{center}
\includegraphics[scale=0.7]{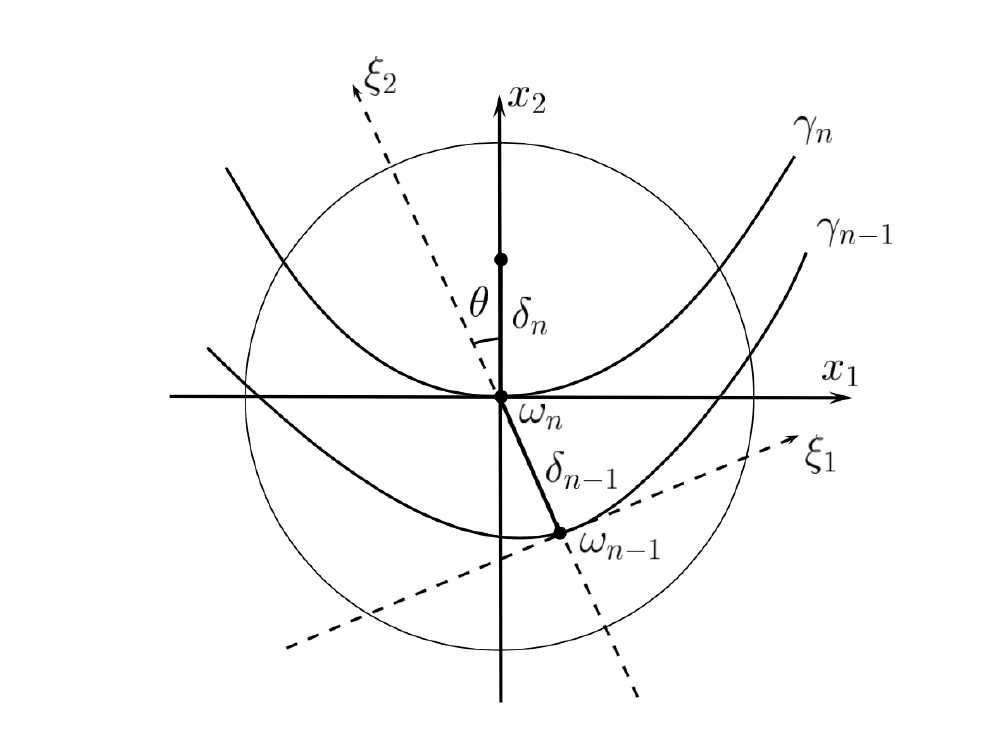}
\caption{The $\xi$-coordinate system.\label{xi-coords}}
\end{center}
\end{figure}

In order to find a precise expansion of the signed distance function $d_{n-1}$ to the previous interface,
it is necessary to calculate the change in the direction of the outer normal to the interface.
Let us denote the point on the previous interface $\gamma_{n-1}$ by $\omega_{n-1}$ and the curvature at that point by $\kappa_{n-1}$.
Analogous notation is also used for the next interface $\gamma_n$.
We are now considering $n=1$, but we still write the symbol $n$ since we will be able to choose $a$ and $b$ to propagate the interface's velocity over each time step. 
According to (\ref{sol1step}), the solution of the wave equation (\ref{wave_eq}) with initial condition $d_{n-1}$ in the $\xi$-coordinate system of figure \ref{xi-coords} reads
\begin{eqnarray*}
u_{n-1}(t,\xi) &=&  \xi_2 + \frac{\kappa}{2} (c^2t^2+\xi_1^2) + \frac{\kappa_{x_1} \xi_1}{6} (3c^2t^2 + \xi_1^2) -\frac{\kappa_{n-1}^2 \xi_2}{2} (c^2t^2 + \xi_1^2) -v_0(0)t - v_{0_{x_1}}(0)t\xi_1 \\
&& \qquad \qquad + \sum_{|\alpha| \geq 3} e_{\alpha}(t,\xi)(t,\xi)^{\alpha} .
\end{eqnarray*}
We want to compute the unit outer normal $\nu_n$ to the level set of $u_{n-1}$ at the point $(0,d_{n-1})$ in the $\xi$-coordinate system.
Hence,
$$ \nu_n = \frac{ \nabla u_{n-1} (\tau,0,\delta_{n-1})}{ | \nabla u_{n-1} (\tau,0,\delta_{n-1})| }, $$
where
\begin{eqnarray*}
\nabla u_{n-1} (\tau,\xi) &=& \left( \kappa_{n-1}\xi_1 + \frac{(\kappa_{n-1})_{x_1}}{2} (c^2 \tau^2+\xi_1^2) - \kappa_{n-1}^2 \xi_1 \xi_2 - v_{0_{x_1}}(0)t,
1- \frac{\kappa_{n-1}^2}{2} (c^2 \tau^2+ \xi_1^2)  \right)^T \\
&& \qquad \qquad \qquad + \nabla_{\xi} \sum_{|\alpha| \geq 3} e_{\alpha}(t,\xi)(t,\xi)^{\alpha} . 
\end{eqnarray*}

This yields
\begin{eqnarray}
\nabla u_{n-1}(\tau,0,\delta_{n-1}) &=& \left( \frac{(\kappa_{n-1})_{x_1}}{2} c^2 \tau^2 - v_{0_{x_1}}(0) \tau, 
1- \frac{\kappa_{n-1}^2}{2} c^2 \tau^2 \right)^T \nonumber \\
\label{normn} 
&& \qquad \qquad 
+ \Big( \sum_{|\alpha| \geq 2} e_{\alpha}(\tau,\delta_{n-1})(\tau,\delta_{n-1})^{\alpha} , \sum_{|\alpha| \geq 2} e_{\alpha}(\tau,\delta_{n-1})(\tau,\delta_{n-1})^{\alpha} \Big)^T,
\end{eqnarray}
and
\begin{equation}
\label{nudef}
\nu_n = \frac{1}{\sqrt{1+(\nu_n^1)^2}}(\nu_n^1,1)^T, \qquad \text{where} \;\;\; 
\nu_n^1 = -v_{0_{x_1}}(0) \tau + \sum_{|\alpha| \geq 2} e_{\alpha}(\tau,\delta_{n-1})(\tau,\delta_{n-1})^{\alpha} . 
\end{equation}

From (\ref{normn}) we can calculate the deviation angle of the normal in one time step,
\begin{eqnarray*} 
\cos \theta = (0,1) \cdot \nu_{n} = 
\frac{1+ \sum_{|\alpha| \geq 2} e_{\alpha}(\tau, \delta_{n-1})(\tau,\delta_{n-1})^{\alpha}}{\sqrt{1+ \sum_{|\alpha| \geq 2} e_{\alpha}(\tau, \delta_{n-1})(\tau,\delta_{n-1})^{\alpha}}} 
= 1+ O(\tau^2) ,
\end{eqnarray*}
since $\delta_{n-1}=O(\tau)$.
Hence, we conclude that $\theta = O(\tau)$.

The $x$-coordinate and $\xi$-coordinate are rotated by the vector $\nu_n$ and thus,
\begin{eqnarray*}
\xi_1 &=& \frac{1}{\sqrt{1+(\nu_n^1)^2}}(x_1+ \nu_n^1 x_2) \\
\xi_2 &=& \frac{1}{\sqrt{1+(\nu_n^1)^2}}(- \nu_n^1 x_1+ x_2) + \delta_{n-1} .
\end{eqnarray*}
We can now write the signed distance $d_{n-1}$ in terms of $x$, as follows (for simplicity we write $\nu$, instead of $\nu_n^1$):
\begin{eqnarray*}
d_{n-1}(x) &=& \xi_2 + \frac{1}{2} \kappa_{n-1} \xi_1^2 +  \frac{1}{6} (\kappa_{n-1})_{x_1}\xi_1^3 - \frac{1}{2} \kappa_{n-1}^2 \xi_1^2 \xi_2 + 
\sum_{| \alpha | = 4} e_{\alpha}(\xi) \xi^{\alpha} \\
&=& \tfrac{1}{\sqrt{1+\nu^2}} (-\nu x_1+x_2) + \delta_{n-1} + \tfrac{\kappa_{n-1}}{2(1+\nu^2)} (x_1+\nu x_2)^2 +  \tfrac{(\kappa_{n-1})_{x_1}}{6\sqrt{1+\nu^2}^3}(x_1+\nu x_2)^3 \\
&& - \tfrac{\kappa_{n-1}^2}{2(1+\nu^2)} (x_1+\nu x_2)^2 (\tfrac{-\nu x_1+x_2}{\sqrt{1+\nu^2}} + \delta_{n-1}) + 
\sum_{| \alpha | = 4} e_{\alpha}(x, \delta_{n-1}) (x,\delta_{n-1})^{\alpha} \\
&=& x_2 + \delta_{n-1} + \frac{1}{2} \kappa_{n-1} x_1^2 +  \frac{1}{6} (\kappa_{n-1})_{x_1}x_1^3 - \frac{1}{2} \kappa_{n-1}^2 x_1^2 x_2 \\
&& - \nu x_1 - \tfrac{\nu^2}{2}x_2 + \tfrac{\kappa_{n-1}}{2} (-\nu^2x_1^2+2\nu x_1 x_2 + \nu^2  x_2^2) 
 + \tfrac{(\kappa_{n-1})_{x_1}}{6} (3\nu x_1^2x_2+3\nu^2 x_1 x_2^2)  
 - \tfrac{(\kappa_{n-1})_{x_1}}{4} \nu^2 x_1^3  \\
&& -\tfrac{\kappa_{n-1}^2}{2} (-\nu x_1^3 -2 \nu^2x_1^2x_2 +2 \nu x_1x_2^2 +\nu^2 x_2^3)
+\tfrac{3\kappa_{n-1}^2}{4} \nu^2 x_1^2x_2 
 -\tfrac{\kappa_{n-1}^2}{2(1+\nu^2)}(x_1+\nu x_2)^2 \delta_{n-1} \\
&& + \sum_{| \alpha | = 4} e_{\alpha}(x, \delta_{n-1}) (x,\delta_{n-1})^{\alpha} + \sum_{i \geq 3} e_i(x) \nu^i .
\end{eqnarray*}
Because of (\ref{nudef}), the last error term can be written as
$$ \sum_{|\alpha| \geq 3} e_{\alpha}(\tau, \delta_{n-1},x)(\tau,\delta_{n-1})^{\alpha} . $$

Computing the solution of the wave equation with initial condition $d_{n-1}$ and zero initial velocity, we obtain
\begin{eqnarray*}
u_{n-1}(t,x) &=&  x_2 + \delta_{n-1}+ \frac{\kappa_{n-1}}{2} (c^2t^2+x_1^2) + \frac{(\kappa_{n-1})_{x_1} x_1}{6} (3c^2t^2 + x_1^2) -\frac{\kappa_{n-1}^2x_2}{2} (c^2t^2 + x_1^2)  \\
&& - \nu x_1 -  \tfrac{\nu^2}{2}x_2 + \tfrac{\kappa_{n-1}}{2} [-\nu^2(c^2t^2+x_1^2)+2\nu x_1 x_2 + \nu^2(c^2t^2+  x_2^2)] \\
&& + \tfrac{(\kappa_{n-1})_{x_1}}{6} [3\nu x_2(c^2t^2+x_1^2)+3\nu^2 x_1(c^2t^2+x_2^2)] - 
\tfrac{(\kappa_{n-1})_{x_1}}{4} \nu^2 x_1(3c^2t^2+x_1^2) \\
&& -\tfrac{\kappa_{n-1}^2}{2} [-\nu x_1(3c^2t^2+x_1^2) -\tfrac{7}{2} \nu^2x_2(c^2t^2+x_1^2) +2 \nu x_1(c^2t^2+x_2^2) +\nu^2 x_2(3c^2t^2+x_2^2)] \\
&& -\tfrac{\kappa_{n-1}^2}{2(1+\nu^2)} \delta_{n-1} [c^2t^2+x_1^2+2\nu x_1x_2 + \nu^2(c^2t^2+x_2^2)] \\
&& + \sum_{|\alpha| =4} e_{\alpha}(t,x)(t,x)^{\alpha} + 
\sum_{|\alpha| \geq 3} e_{\alpha}(\tau, \delta_{n-1},t,x)(\tau,\delta_{n-1})^{\alpha},
\end{eqnarray*}
where we have used the calculations of the integrals from Poisson formula in Section \ref{sec1BMO}.
The solution for for the initial condition $ad_{n-1}+bd_n$ evaluated at the point $(t,x)=(\tau,0,\delta_n)$ thus reads
\begin{eqnarray*}
u_{n}(\tau,0,\delta_n) &=&  (a+b)\delta_{n} + a \delta_{n-1} + \frac{a\kappa_{n-1}+b\kappa_n}{2} c^2\tau^2 
-\frac{a \kappa_{n-1}^2 + b \kappa_n^2}{2} \delta_n c^2\tau^2 \\
&& + a \Big\{ - \tfrac{\nu^2}{2} \delta_n + \tfrac{\kappa_{n-1}}{2} \nu^2 \delta_n^2 
+ \tfrac{(\kappa_{n-1})_{x_1}}{2} \nu \delta_n c^2\tau^2 \\
&& -\tfrac{\kappa_{n-1}^2}{2} [-\tfrac{1}{2} \nu^2 \delta_{n} c^2\tau^2 + \nu^2 \delta_{n}^3] 
-\tfrac{\kappa_{n-1}^2}{2(1+\nu^2)} \delta_{n-1} [c^2\tau^2 (1+\nu^2) + \nu^2 \delta_{n}^2] \Big\} \\
&& + \sum_{|\alpha| =4} e_{\alpha}(\tau,\delta_n)(\tau,\delta_n)^{\alpha} + 
\sum_{|\alpha| \geq 3} e_{\alpha}(\tau, \delta_{n-1},\delta_n)(\tau,\delta_{n-1})^{\alpha} .
\end{eqnarray*}

The distance $\delta_{n}$ traveled by the interface in the normal direction then satisfies $u_n(t,0,\delta_n)=0$, which can be written as
\begin{eqnarray*}
0 &=& - \delta_{n}^3 a \nu^2 \tfrac{\kappa_{n-1}^2}{2}  + \delta_{n}^2 a \nu^2 \tfrac{\kappa_{n-1}^2}{2} (1- \tfrac{\delta_{n-1}}{1+\nu^2}) \\
&& + \delta_{n} (a+b - \tfrac{a\kappa_{n-1}^2+b\kappa_n^2}{2} c^2\tau^2 - \tfrac{a}{2} \nu^2 + \tfrac{a}{2} \nu c^2\tau^2 ((\kappa_{n-1})_{x_1}  + \tfrac{1}{2} \nu \kappa_{n-1}^2)) \\
&& + a \delta_{n-1} + \tfrac{a\kappa_{n-1}+b\kappa_n}{2} c^2\tau^2 - \tfrac{\kappa_{n-1}^2}{2} a \delta_{n-1} c^2\tau^2 + \text{error terms} ,
\end{eqnarray*}
where the error terms are of order $O(\tau^3+\delta_{n-1}^3+\delta_n^4)$.
Moreover, the quantities $\delta_{n-1}$ and $\nu$ are of order $O(\tau)$, at worst. 
Since the lowest order approximation of the above equation is 
\begin{eqnarray*}
\delta_{n} (a+b) + a \delta_{n-1} + \tfrac{a\kappa_{n-1}+b\kappa_n}{2} c^2 \tau^2 =0 ,
\end{eqnarray*}
we see that $\delta_n$ must be of order $O(\delta_{n-1}+\tau^2)$, which is at most $O(\tau)$.
Hence the $\delta_n^2$ and $\delta_n^3$ terms (including their coefficients) are of order $O(\tau^4)$ or higher and can be neglected.
In this way we obtain the solution $\delta_n$ to the equation as
\begin{eqnarray*}
\delta_{n} &=& - \frac{a \delta_{n-1} + \tfrac{a\kappa_{n-1}+b\kappa_n}{2} c^2 \tau^2 + O(\tau^3 + \delta_{n-1}^3)}{a+b - \tfrac{a\kappa_{n-1}^2+b\kappa_n^2}{2} c^2 \tau^2 - \tfrac{a}{2} \nu^2} \\
&=& \tfrac{-a}{a+b} \delta_{n-1} - \tfrac{a\kappa_{n-1}+b\kappa_n}{2(a+b)} c^2 \tau^2 + O(\tau^3) .
\end{eqnarray*}

It is natural to choose $a,b$ so that $-\frac{a}{a+b}=1$ and $a+b=1$, i.e., $a=-1,b=2$.
Moreover, setting $c^2=2$ yields
\begin{equation}\label{acc-curve}
\delta_{n} = \delta_{n-1} - \left( 2 \kappa_{n}- \kappa_{n-1} \right) \tau^2 + O(\tau^3),
\end{equation}
which is the desired relation between $\delta_n$ and $\delta_{n-1}$. In particular, if we assume that the curvature
is changing smoothly, relation (\ref{acc-curve}) can be rewritten as follows,
\begin{eqnarray}
\delta_{n} &=& \delta_{n-1} - \kappa_n \tau^2 - \left( \kappa_{n}- \kappa_{n-1} \right) \tau^2 + O(\tau^3) 
\nonumber \\
\label{obsc}
&=& \delta_{n-1} - \kappa_n \tau^2 + O(\tau^3) .
\end{eqnarray}

To see that this is a correct approximation, consider the one-dimensional point-mass motion with initial velocity $v_0$ and acceleration $-\kappa (t)$,
thus solving $x''(t)= -\kappa (t), \; x(0)=0, \; x'(0)=v_0$. The solution is
\begin{equation}
\label{1dsol}
x(t) = v_0 t - \int_0^t \int_0^s \kappa(u) \, du \, ds . 
\end{equation}
Our algorithm (\ref{obsc}), including the initial step from the previous subsection, thus gives the approximation
$$ \delta_1 = v_0 \tau - \frac{1}{2} \kappa_0 \tau^2 + O(\tau^3) , \; \delta_2 = v_0 \tau - \frac{1}{2} \kappa_0 \tau^2 - \kappa_1 \tau^2 + O(2\tau^3) ,
 \dots , \delta_k = v_0 \tau - \frac{1}{2} \kappa_0 \tau^2 - \tau^2 \sum_{i=1}^{k-1} \kappa_i + O(k \tau^3)  . $$
The corresponding total distance is
\begin{equation}
\label{1djust}
\tilde{x}(k \tau) = kv_0 \tau - k \frac{1}{2} \kappa_0 \tau^2 - \tau^2 \sum_{i=1}^{k-1} \sum_{j=1}^{k-i} \kappa_j + (k \tau)^2 O(\tau) , 
\end{equation}
which converges as $\tau \to 0$ to the function (\ref{1dsol}) by the trapezoidal rule.
The above analysis also shows that the accumulation of errors through time does not spoil the approximation of the interface, at least until the development of singularities.

\noindent{\bf{Remark}}. Regarding the 3D problem, assuming that the signed distance can be expanded as
$$ d(y) = y_3 + \frac{1}{2} \kappa_1 y_1^2 + \frac{1}{2} \kappa_2 y_2^2 + H.O.T., $$
then a similar approach using the solution to the three-dimensional wave equation suggests the validity of our method in higher dimensions:
\begin{eqnarray*}
u(t,x) &=& \frac{1}{4 \pi t^2} \int_{\partial B(x,t)} \Big( d(y) + \nabla d(y) \cdot (y-x) \Big) \, dS(y) \\
&=& \frac{1}{4 \pi t^2} \int_{\partial B(x,t)} \Big( y_3 + \frac{1}{2} \kappa_1 y_1^2 + \frac{1}{2} \kappa_2 y_2^2 + \kappa_1y_1(y_1-x_1)+
\kappa_2 y_2 (y_2-x_2) + (y_3-x_3) \Big) \, dS(y) \\
&=& \frac{1}{4 \pi} \int_{\partial B(0,1)} \Big( tz_3 + x_3 + \frac{1}{2} \kappa_1 (tz_1+x_1)^2 + \frac{1}{2} \kappa_2 (tz_2+x_2)^2 \\
&& \qquad \qquad \qquad + \kappa_1tz_1(tz_1+x_1)+ \kappa_2 tz_2 (tz_2+x_2) + tz_3 \Big) \, dS(z) \\
&=& \frac{1}{4 \pi} \int_{\partial B(0,1)} \Big( x_3 + \frac{1}{2} \kappa_1 (3t^2z_1^2+x_1^2) + \frac{1}{2} \kappa_2 (3t^2z_2^2+x_2^2) \Big) \, dS(z) \\
&=& \frac{1}{4 \pi} \int_0^{2 \pi} \int_0^{\pi} \Big( x_3 + \frac{1}{2} \kappa_1 x_1^2 + \frac{1}{2} \kappa_2 x_2^2 + \frac{3}{2} \kappa_1 t^2 \sin^2 \theta \cos^2 \phi + \frac{3}{2} \kappa_2 t^2 \sin^2 \theta \sin^2 \phi \Big) \sin \theta \, d \theta \, d \phi \\
&=& \frac{1}{4 \pi} \left[ 2 \cdot 2 \pi \big( x_3 + \frac{1}{2} \kappa_1 x_1^2 + \frac{1}{2} \kappa_2 x_2^2 \big) + \frac{3}{2} (\kappa_1 + \kappa_2) t^2 \frac{4}{3} \cdot \pi \right]  \\
&=& x_3 + \frac{1}{2} \kappa_1 (t^2+x_1^2) + \frac{1}{2} \kappa_2 (t^2+x_2^2).
\end{eqnarray*}
This yields the distance traveled by interface in time $t$ as
$$ \alpha = - \frac{1}{2} (\kappa_1 + \kappa_2) t^2 = - H t^2 ,$$
confirming that our scheme can also be used with the three-dimensional wave equation, and likely in any higher dimension.
\section{Numerical Tests and Properties}

Thresholding dynamical algorithms have the computational advantage that there is no need to calculate curvatures, and singularities are implicitly handled by the partial differential equation. Moreover, they are extremely simple to implement (here, we need only solve the wave equation). Thus, one can construct the corresponding numerical method in a number of ways--we will make use of standard finite differences, as well as {\emph{minimizing movements}} for use in investigating volume preserving motions.

Our approximation method for computing the HMCF (\ref{hmcf}) requires one to solve a wave equation, for which there are a number of schemes. Since we are considering an interfacial motion embedded in the solution to a wave equation, care needs to be taken when detecting the interface and constructing the signed distance functions. In particular, the precise location of the interface is needed. Without this information, errors arising from its  approximation will propagate via the evolution of the signed distance functions.

In the same vein,  since the governing partial differential equation is hyperbolic, 
care must be taken in the numerical implementation of the method when constructing initial conditions. In particular, phase locations are encoded by 
vectors at the nodes of the discretized domain, and so the finite element assumptions can introduce initial interfaces that are not smooth. One approach to alleviating issues caused by this is to use one step of the original BMO algorithm to smooth the initial interface. This is easily done, since it only requires one to solve the vector valued heat equation for a small time.

We will first perform a convergence test for an idealized version of our algorithm. This test assumes that no errors are introduced under spacial discretization. Our second test introduces the spacial discretization, but computes signed distance functions utilizing an idealized representation of the interface. The final test removes this idealization and examines our algorithm's order of convergence using standard finite differences. 

\subsection{Idealized numerical convergence}
This section investigates the convergence rate of our algorithm for a simple test problem, under a slight idealization.
For a circle evolving by (\ref{hmcf}) with initial normal velocity $v_0$ (we assume $v_0 \leq 0$ for simplicity), one can compute the evolution of the circle's radius $r(t)$ analytically by solving
\begin{equation}\label{ode}
 r''(t)= - \frac{1}{r(t)} , \qquad r(0) = r_0, \qquad r'(0) = v_0.
 \end{equation}
Multiplying the equation by $r'(t)$ and integrating we obtain 
\begin{eqnarray*}
r'(t) &=& - \sqrt{-2 \log r + C_1}, 
\end{eqnarray*}
where $C_1 = 2 \log r_0 + v_0^2$ is determined by the initial conditions.
Since
$$ \int \frac{dr}{\sqrt{-2 \log r +C_1}} = - e^{C_1/2} \sqrt{\tfrac{\pi}{2}} \text{erf} \left( \sqrt{ \tfrac{C_1}{2} - \log r} \right) , $$
we have
$$ e^{C_1/2} \sqrt{\tfrac{\pi}{2}} \text{erf} \left( \sqrt{ \tfrac{C_1}{2} - \log r(t)} \right) = t + C_2 , \qquad
\text{where} \;\; C_2 = r_0 \sqrt{\tfrac{\pi}{2}} e^{\frac{v_0^2}{2}} \text{erf} \Big( \frac{v_0}{\sqrt{2}} \Big) . $$
In particular, $C_2=0$ for zero intial velocity. 

The solution can thus be expressed as
$$ r(t) = \text{exp} \left( \frac{C_1}{2} - \left[ \text{erf}^{-1} \Big( \sqrt{\tfrac{2}{\pi}} e^{-C_1/2} (t+C_2) \Big) \right]^2 \right) , $$
which yields a more simple form for zero initial velocity:
$$ r(t) = r_0 e^{- \left[ \text{erf}^{-1} \left( \sqrt{\frac{2}{\pi}} \frac{t}{r_0} \right) \right]^2 } . $$

\begin{figure}[h!]
\begin{center}
\includegraphics[scale=0.5]{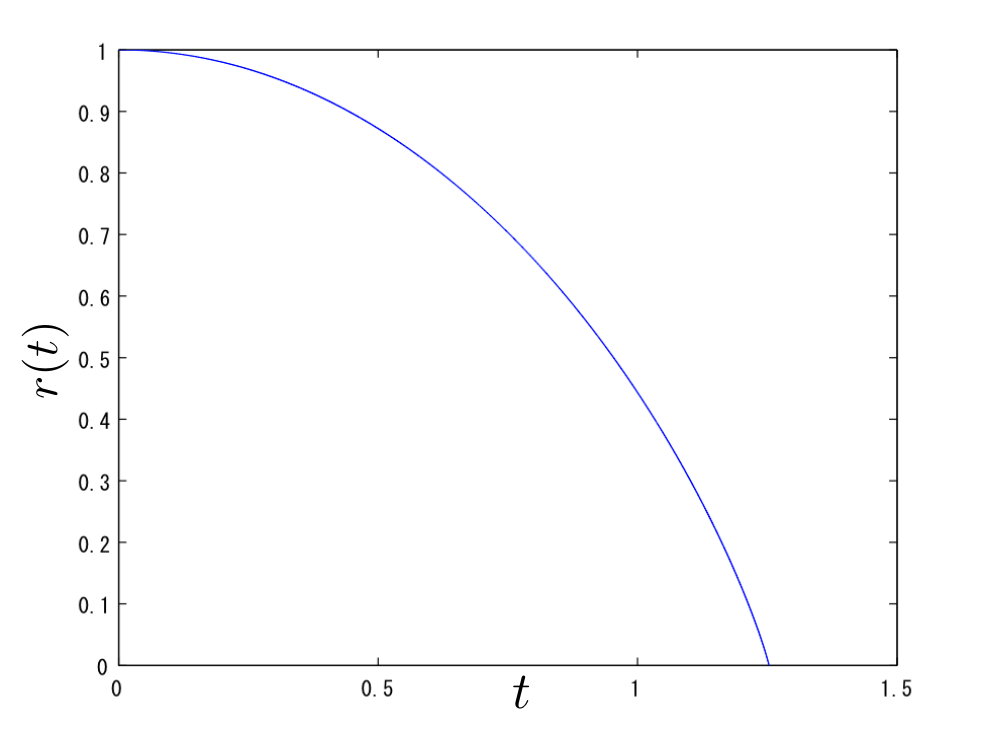}
\caption{Graph of the solution for $r_0=1, v_0=0$.\label{hmcfode}}
\end{center}
\end{figure}

Since the function $\text{erf}^{-1}(x)$ tends to infinity when $x \to 1-$, the radius vanishes in a finite time, which can be 
computed from the relation
$$ \sqrt{\tfrac{2}{\pi}} e^{-C_1/2} (t_e+C_2) = 1. $$
After substituting values of integration constants, this gives
\begin{equation}\label{extinction}
 	t_e = r_0 \sqrt{\tfrac{\pi}{2}} e^{v_0^2/2} \Big( 1+ \text{erf} \big( \frac{v_0}{\sqrt{2}} \big) \Big) .
\end{equation}
We remark that the extinction time for zero initial velocity is $t_e = r_0 \sqrt{\pi /2}$ (see figure \ref{hmcfode}).

Let us now calculate the numerical approximation due to our scheme.
That is, we consider two circles of radii $r_0 > r_1$, centered at the origin, which describe the initial conditions. 
We then compute the radius of the circle given by the $0$-level set of the solution (after time $t$) to the problem 
$$ u_{tt} = c^2 \Delta u, \qquad u(0,x)=2d_1(x)-d_0(x), \qquad u_t(0,x)=0 . $$
Here, $c^2=2$ and $d_i$ is the signed distance function to the circle of radius $r_i$, which can be written as $d_i(x)=r_i - |x|$.
Hence, denoting $r=2r_1-r_0$, the initial condition is given by $u(0,x)=r-|x|$ 
and the solution of the above initial value problem (sufficiently far away from the origin) is
\begin{eqnarray*}
u(t,x) = \frac{1}{2c\pi t} \int_{B(x,ct)} \frac{r-|y|- \frac{y}{|y|} \cdot (y-x)}{\sqrt{c^2t^2-|y-x|^2}} \, dy \
= \frac{1}{2\pi} \int_{B(0,1)} \frac{r-2|x+ctz| + \frac{|x|^2 + ct x \cdot z}{|x+ ctz|}}{\sqrt{1-|z|^2}} \, dz .
\end{eqnarray*}
Since the initial condition is radially symmetric and there is no interaction with a boundary, we can restrict the values of
$x$ to points on the $x$-axis of the form
$x=(\xi,0)$, where $\xi> 0$.
Then the above simplifies to
\begin{eqnarray*}
u(t,\xi,0) 
= r - \frac{1}{2\pi} \int_0^1 \int_0^{2\pi} \frac{\xi^2+3ct \rho \xi \cos \theta+2c^2 t^2 \rho^2}{\sqrt{1-\rho^2} \sqrt{\xi^2+2ct \rho \xi \cos \theta + c^2 t^2 \rho^2}} \rho \, d \theta \, d \rho .
\end{eqnarray*}
The exact form of the solution to the above equation is difficult to find, since it involves elliptic integrals.
Therefore, let us instead assume that $t=\tau$ is small and expand the integrand as a Taylor series in time.
To this end, let us denote
\begin{eqnarray*}
a(t) = \xi^2+3ct \rho \cos \theta+2c^2 t^2 \rho^2, \quad
b(t) = \sqrt{\xi^2+2ct \rho \xi \cos \theta + c^2 t^2 \rho^2} , \quad
G(t) = \frac{a(t)}{b(t)}, 
\end{eqnarray*}
and calculate the first and second derivatives of $G$ at $t=0$.
\begin{eqnarray*}
G(0) &=& \xi , \\
G'(t) 
&=&  \frac{2c \xi^3 \rho \cos \theta + 3c^2 \xi^2 t \rho^2 (1+ \cos^2 \theta) + 8c^3 \xi t^2 \rho^3 \cos \theta + 2 c^4 t^3 \rho^4}{ b^3(t)} , \\
G'(0) &=& 2c \rho \cos \theta , \\
G''(0) &=& \frac{3c^2 \xi^2 \rho^2 (1+\cos^2 \theta) \cdot \xi^3 - 2c \xi^3 \rho \cos \theta \cdot \frac{3}{2} \xi 2c \rho \xi \cos \theta}{\xi^6} 
= \frac{3c^2}{\xi} \rho^2 \sin^2 \theta .
\end{eqnarray*}
We can thus approximate the solution at points $(\xi,0)$ as
\begin{eqnarray*}
u(\tau,\xi,0) &=& r - \frac{1}{2 \pi} \int_0^1 \int_0^{2 \pi} \frac{\rho}{\sqrt{1- \rho^2}} \Big( G(0) + G'(0)\tau + \tfrac{1}{2} G''(0)\tau^2 +O(\tau^3) \Big) \, d \theta \, d \rho \\
&=& r - \frac{1}{2 \pi} \int_0^1 \int_0^{2 \pi} \frac{\rho}{\sqrt{1- \rho^2}} \Big( \xi + 2 c \tau \rho \cos \theta + \tfrac{3c^2}{2\xi} \tau^2 \rho^2 \sin^2 \theta +O(\tau^3) \Big) \, d \theta \, d \rho \\
&=& r -\xi - \frac{c^2}{2\xi} \tau^2 + O(\tau^3).
\end{eqnarray*}
To determine the position $\xi$ of the interface after time $\tau$, we need to solve the equation
$$ r -\xi - \frac{c^2}{2\xi} \tau^2 = O(\tau^3) . $$
The corresponding quadratic equation $ \xi^2-(r+O(\tau^3))\xi + \tfrac{1}{2} c^2 \tau^2 = 0$ has two roots 
but since for vanishing times we want to recover the value $\xi=r$, we take the corresponding root
\begin{eqnarray*}
\xi &=& \frac{1}{2} \left( r + \sqrt{r^2 -2c^2 \tau^2} \right) + O(\tau^3) 
= \frac{1}{2} \left( r + r - \frac{c^2 \tau^2}{r} \right) + O(\tau^3) 
= r - \frac{1}{r} \tau^2 + O(\tau^3) .
\end{eqnarray*} 

Combining this result with our considerations regarding the first HBMO step, we obtain the following HBMO algorithm for the special case of a circle.
\begin{enumerate}
\item Prescribe initial radius $r_0$, initial velocity $v_0$ and time step $\tau=t_e/N$, where $t_e$ is the extinction time (\ref{extinction}) and $N \in {\mathbb N}$.
\item Set 
$$ r_1 = r_0- \frac{\tau^2}{2r_0} . $$
\item Repeat the following for $n=2, \dots ,N$:
\begin{equation}
\label{hbmo_circ}
r_n = 2r_{n-1} - r_{n-2} - \frac{\tau^2}{2r_{n-1}-r_{n-2}} . 
\end{equation}
\end{enumerate}  

We test the resulting scheme for $r_0=1$ and $v_0=0$.
The results are shown in the table below, where the error is taken at the half time to extinction, i.e.,
$$ \text{error} = |r_{N/2} - r(t_e/2)|=|r_{N/2} - exp(-[erf^{-1}(\tfrac{1}{2})]^2)|. $$
We observe exactly linear convergence with respect to the time step $\tau$ as predicted by (\ref{1djust}), 
which shows that the proposed HBMO algorithm gives correct (i.e., convergent) result when the initial curve is a circle.

\begin{table}[h!]
\begin{center}
  \begin{tabular}{|l|l|l|} \hline
   division number $N$ & error at time $t_e/2$ & convergence order \\ \hline \hline
		10     & 0.073199 &  - \\ \hline
		$10 \cdot 4^1$ & 0.019332 &  0.960 \\ \hline
		$10 \cdot 4^2$ & 0.004884 &  0.992 \\ \hline
		$10 \cdot 4^3$ & 0.001224 &  0.998 \\ \hline
		$10 \cdot 4^4$ & 0.000306 &  1.000 \\ \hline
		$10 \cdot 4^5$ & 0.000076 &  1.000 \\ \hline
		$10 \cdot 4^6$ & 0.000019 &  1.000 \\ \hline
		$10 \cdot 4^7$ & 0.000004 &  1.000 \\ \hline
	\end{tabular}
	\caption{Results of the HBMO algorithm (\ref{hbmo_circ}) with discretized time only.}
\end{center}
\end{table}

\subsection{Convergence analysis under spacial discretization}

Our next tests concern the algorithm's convergence under finite differencing. A circle of radius $3/10$ is centered at $(x,y)=(1/2,1/2)$ with zero initial velocity, and the initial condition for our algorithm is therefore the signed distance function (taking positive values inside the circle) to this interface. The calculations are performed in the unit square under zero Neumann boundary conditions and we use grid points on a square lattice with spacing $\Delta x = \Delta y = 1/(N-1)$. A Delaunay triangulation of these points is performed, and we assume a linear interpolation of the grid values within each element. This allows us to define and determine the precise location of the interface. In particular, the location of the interface at each time step is given as a collection of line segments, each of which corresponds to the zero level set of the signed distance function across the elements. We remark that the geometry of the interface depends on the triangulation and that, as the interface shrinks, the relative resolution of the grid is a source of error. 

We take the HBMO time step as $\Delta t = t_e/{2^9}$ ($t_e \approx 0.376$) and the time step for the finite difference approximation of the solution to the
wave equation is $\tau = \Delta t/(2^6)$. We then compute as follows:

\begin{enumerate}
\item Construct the signed distance function $d_0(x,y)$ to the initial circle.
\item Set $d_1(x,y)=d_0(x,y)$ (zero initial velocity), and repeat steps $(3-5)$, for $n=1,\dots,MAX$.
\item Compute the finite difference approximation to the following the wave equation:
\begin{equation}\label{discretewave1}
\begin{cases}
u_{tt} = \Delta u \hspace{119pt} \text{in } (0,\Delta t)\times\Omega\\
\frac{\partial u}{\partial \nu}=0 \hspace{121pt} \text{on } (0,\Delta t)\times\partial\Omega\\
u(t=0,x) = 2d_{n}-d_{n-1},\hspace{10pt} u_t(t=0,x)=0 \hspace{11pt} \text{in } \Omega.
\end{cases}
\end{equation}
\item Compute $\bar{r}_{n+1}$, the average distance from $(1/2,1/2)$ of the line segments describing the zero level set of the solution to (\ref{discretewave1}).
\item Set $d_{n+1}(x,y) =  \sqrt{(1/2 - x)^2 + (1/2 - y)^2}-\bar{r}_{n+1}.$
\end{enumerate}  

The error table is show in Table \ref{errortab_fin_diff_ideal}. Here, $L^2$-error designates to the square root of the value below, 
\begin{equation}
\Delta t \sum_{n=1}^{N_{e}} \left( r((n-1)t) - \bar{r}_{n}\right)^2,\notag
\end{equation}
where $r(t)$ denotes the solution to (\ref{ode}), $N_{e}$ is the largest positive integer satisfying $N_{e}\Delta t \leq t_e$, and the value of $\bar{r}_{n}$ is zero for any iteration after which the approximate solution's radius becomes zero.

\begin{table}[h!]
\begin{center}
  \begin{tabular}{|l|l|l|} \hline
   Grid Resolution $N$ & $L^2$-error & convergence order \\ \hline \hline
		16        & 0.132320    &  - \\ \hline
		$32$    & 0.126222    & 0.068067\\ \hline
		$64$    & 0.113587    & 0.152170 \\ \hline
		$128$  & 0.088544    & 0.359340  \\ \hline
		$256$  & 0.051558    & 0.780184\\ \hline
		$512$  & 0.022104    & 1.221893\\ \hline
		$1024$& 0.006538    & 1.757309\\ \hline
	\end{tabular}
	\caption{Errors obtained for finite differences (using ideal distance functions).\label{errortab_fin_diff_ideal}}
\end{center}
\end{table} 

As our final numerical test, we replace the steps $4$ and $5$ of the previous algorithm with the construction of the exact signed distance function according to the exact configuration of the finite element space's approximation of the interface.
In both cases, we perform the numerical computations on various grid resolutions and investigate the corresponding convergence orders. 


The error table is show in Table \ref{errortab_fin_diff}, where we use the same convention as in the previous test. We remark that explicit computation of the distance functions does not change the order of convergence.

\begin{figure}[h!]
	\begin{center}
		{\includegraphics[scale=0.35,trim=1cm 0cm 1cm 0cm,clip=]{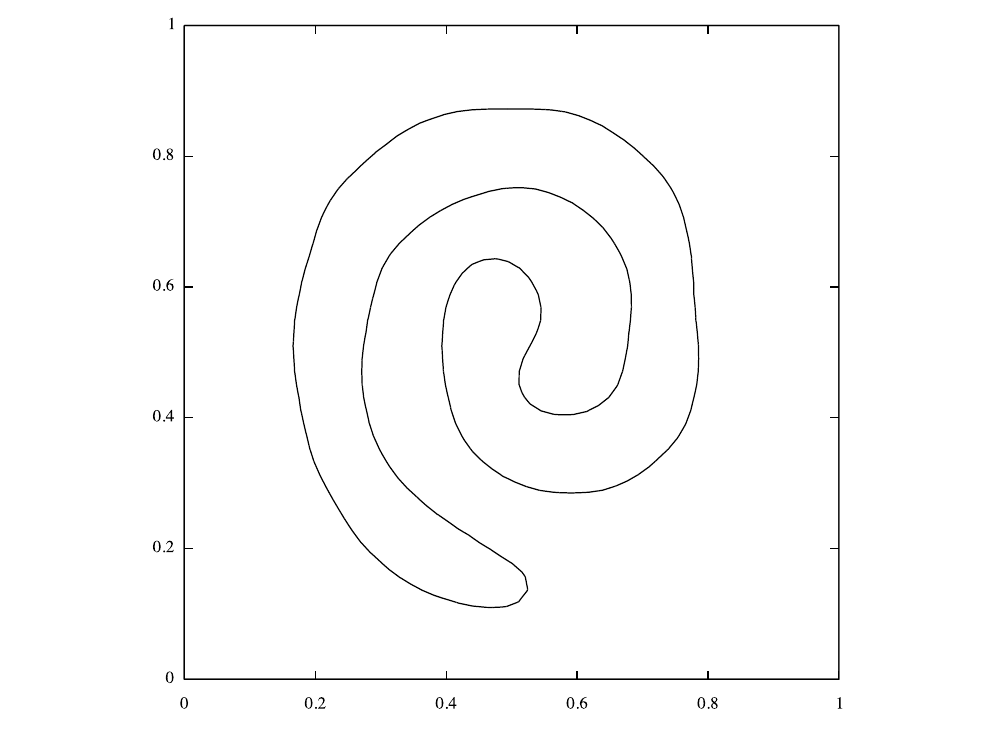}	}
		{\includegraphics[scale=0.35,trim=1cm 0cm 1cm 0cm,clip=]{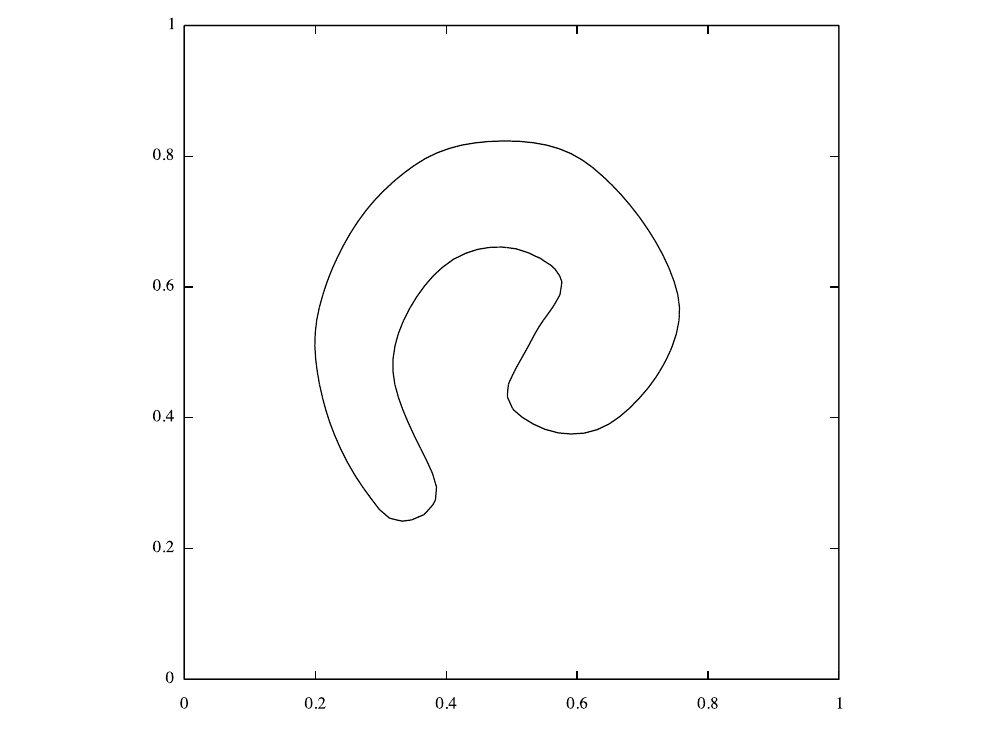}	}
		{\includegraphics[scale=0.35,trim=1cm 0cm 1cm 0cm,clip=]{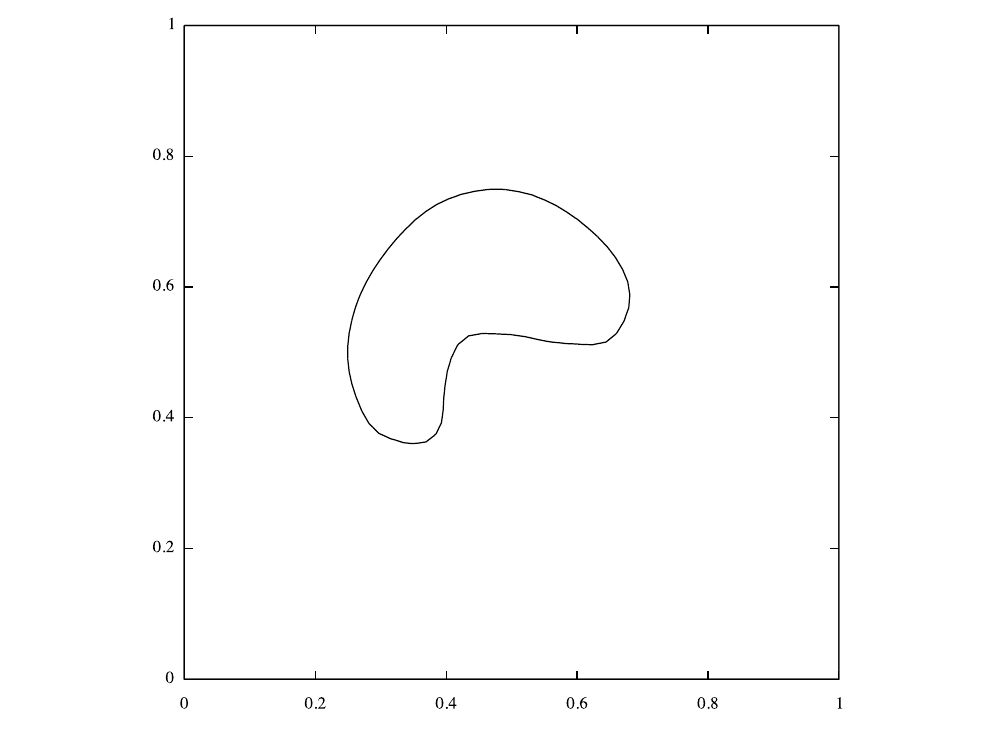}	}
		{\includegraphics[scale=0.35,trim=1cm 0cm 1cm 0cm,clip=]{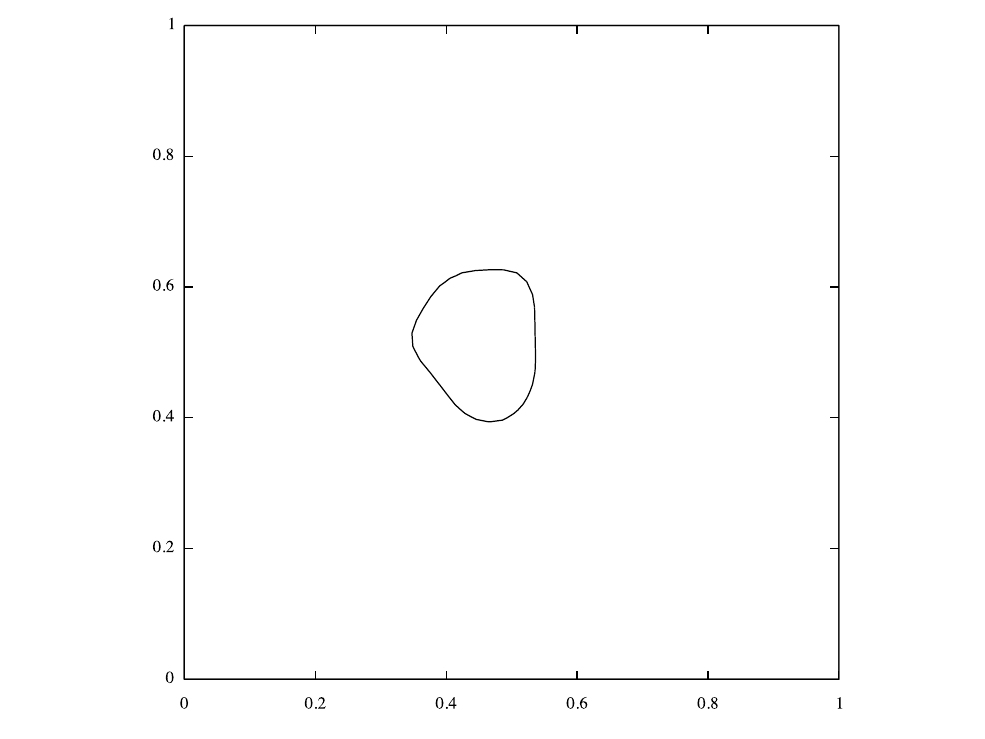}	}
		{\includegraphics[scale=0.35,trim=1cm 0cm 1cm 0cm,clip=]{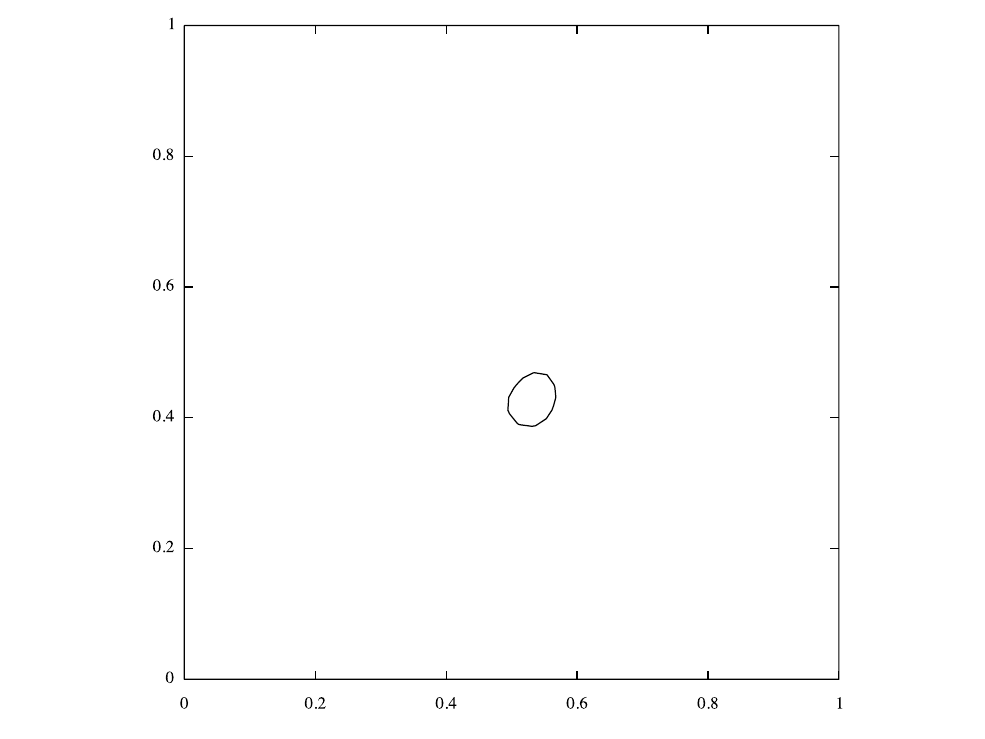}	}
		\caption{Evolution by $\bu'=-\kappa \bnu$ (time is from left to right).}
		\label{hmcffig1}
	 \end{center}
\end{figure}

\begin{table}[h!]
\begin{center}
  \begin{tabular}{|l|l|l|} \hline
   Grid Resolution $N$ & $L^2$-error & convergence order \\ \hline \hline
		16       &   0.131575 &  - \\ \hline
		$32$   &   0.124060 & 0.084851\\ \hline
		$64$   &   0.109666 & 0.177921 \\ \hline
		$128$ &   0.084465 & 0.376696 \\ \hline
		$256$ &   0.044312 & 0.930634 \\ \hline
		$512$ &   0.021216 & 1.062548 \\ \hline
	\end{tabular}
	\caption{Errors obtained for finite differences. \label{errortab_fin_diff}}
\end{center}
\end{table}         
Figure \ref{hmcffig1} shows the evolution of a two phase HMCF, with zero initial velocity, as obtained by our algorithm. Contrast to parabolic mean curvature flow, we observe oscillations in the interfacial motion--the effect of acceleration is also evident in the motion of the phase region's centre. 

\subsection{Multiphase motions}
Since the solution to the wave equation depends only on local information, one can formulate a multiphase algorithm which is akin to the original multiphase BMO algorithm in \cite{[mbo]}. In particular, thresholding a collection of independently evolving wave equations, can formally be shown to compute multiphase interfacial motions with normal acceleration equal to curvature. Nevertheless, since we would also like to investigate multiphase volume preserving motions, we introduce a reformulation of this idea. In particular, by choosing a small time step $\Delta t$, we find a function $\bu:\Omega\rightarrow {\bf{R}}^{N-1}$ solving the following vector valued wave equation:
\begin{equation}\label{vpde}
\begin{cases}
	\bu_{tt} = \Delta \bu \hspace{92pt} \text{ in } (0,\Delta t)\times\Omega,\\
	\frac{\partial \bu}{\partial \nu}=\boldsymbol{0} \hspace{96pt} \text{ on } (0,\Delta t)\times\partial\Omega,\\
	\bu_{t}(0,x)=\boldsymbol{0}\hspace{125pt} \text{ in } \Omega,\\
	\bu(t=0,x) = 2{\zu}^\epsilon_{0} - {\zu}^\epsilon_{-\Delta t} \hspace{67pt} \text{ in } \Omega,
\end{cases}
\end{equation}
where $N$ denotes the number of phases, $\Omega$ is a smooth bounded domain in ${\bf{R}}^d$, 
and the initial condition is defined by the following signed-distance interpolated vector field (see \cite{[eg-thesis]} and \cite{[rhudaina]})
\begin{align}\label{sdinterp}
\zu^{\epsilon}_t(x) =  \sum_{i=1}^{N} \pu_{i}\chi_{\{d_{i}^t(x) > \eps/2\}} + \frac{1}{\eps}\left( \frac{\eps}{2} 
	+ d_{i}^t(x) \right) \pu_{i}\chi_{\{ -\eps/2 \le d_{i}^t(x) \le \eps/2\}}.
\end{align}
Here, $\chi_{E}$ is the characteristic function of the set $E$, the vector $\pu_{i}$ is the $i^{th}$ coordinate vector of a {\emph{regular simplex}} in ${\bf{R}}^{N-1}$ for $i=1,...,N$, and $d_{i}^t(x)$ denotes the signed distance function to the boundary of phase $i$ at location $x$ and time $t$:
 \begin{align}\label{signeddistance}
 	d_{i}^t(x) = 
	\begin{cases}
		\hspace{11pt}\inf_{y \in \partial P_i^t} \vnorm{x-y} \hspace{25pt}\text{if $x\in P_{i}^t$},\\
		-\inf_{y \in \partial P_i^t} \vnorm{x-y} \hspace{25pt}\text{otherwise}.
	\end{cases}
 \end{align}
 We remark that, when $N=2$, equation (\ref{vpde}) is scalar. 
At time $\Delta t$, in a process called {\emph{thresholding}}, each phase region $P_{i}^{\Delta t}$ is evolved as follows:
\begin{align}\label{trunc}
	P_{i}^{\Delta t} = \{ x \in \Omega : \bu (\Delta t,x)\cdot \pu_i \geq \bu (\Delta t,x)\cdot \pu_k, \text{for all } k\in \{1,...,N\} \}.
\end{align}
The vector field ${\zu^\epsilon_{0}}$ is then reconstructed using the boundaries of these sets and the initial condition for the wave equation is updated. The procedure is then repeated. The defining characteristic of the vector field (\ref{sdinterp}) is that, upon choosing a location $x$ positioned within a distance $\epsilon/2$ of an interface, say corresponding to the boundary of phase $k$,  $\zu^{\epsilon}$ satisfies the property:
\begin{align}\label{scalardist}
	\zu^{\epsilon}_t(x)\cdot \pu_k = d^t_{k}(x),
\end{align} 
so that the profile in the direction normal to the interface is given by a signed distance function. This fact allows one to obtain the curvature of the interface by computing the Laplacian of (\ref{scalardist}). Moreover, as shown in \cite{[GOS3]}, this setting introduces a multiple well potential which allows one to express multiphase volume constrained motions. 

\begin{figure}[htbp]
	\centering
   	\begin{tabular}{cccc}
	\hspace{-10pt}
	\setlength\fboxsep{0pt}
	\setlength\fboxrule{0.5pt}
		{\includegraphics[scale=0.43,trim=1cm 0cm 1cm 0cm,clip=]{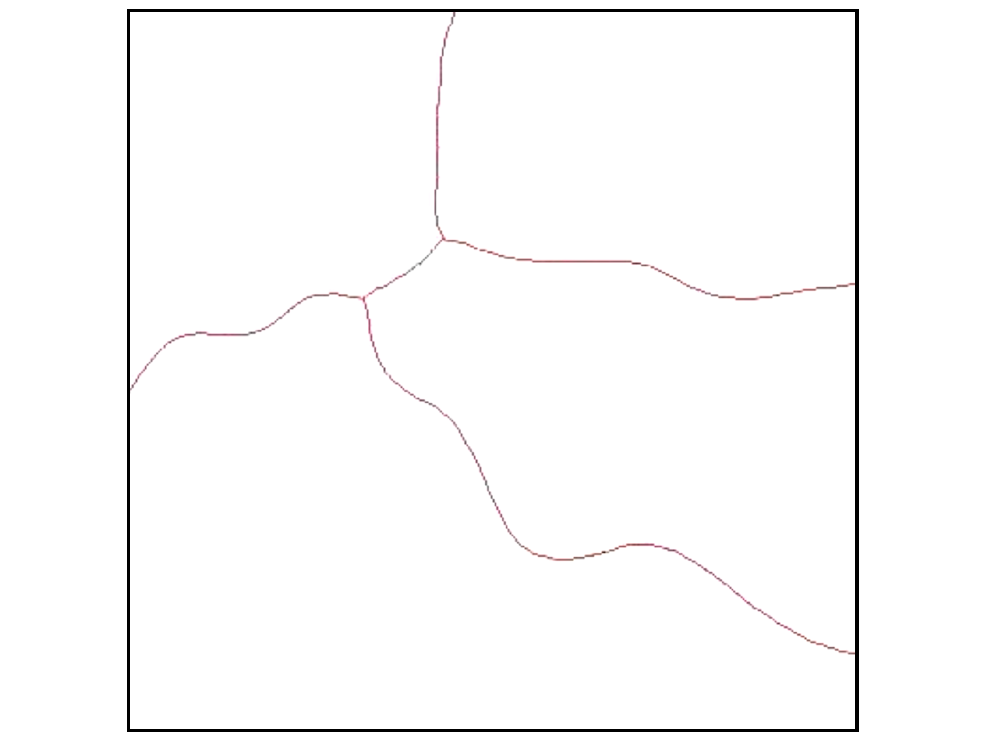} }
		{\includegraphics[scale=0.43,trim=1cm 0cm 1cm 0cm,clip=]{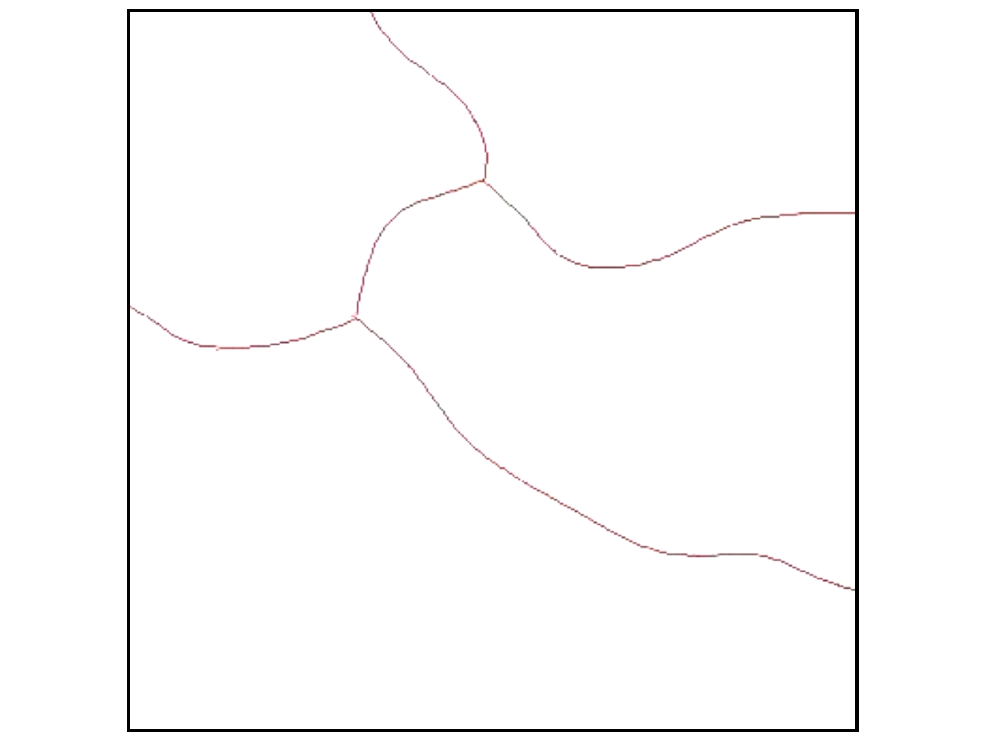} }
		{\includegraphics[scale=0.43,trim=1cm 0cm 1cm 0cm,clip=]{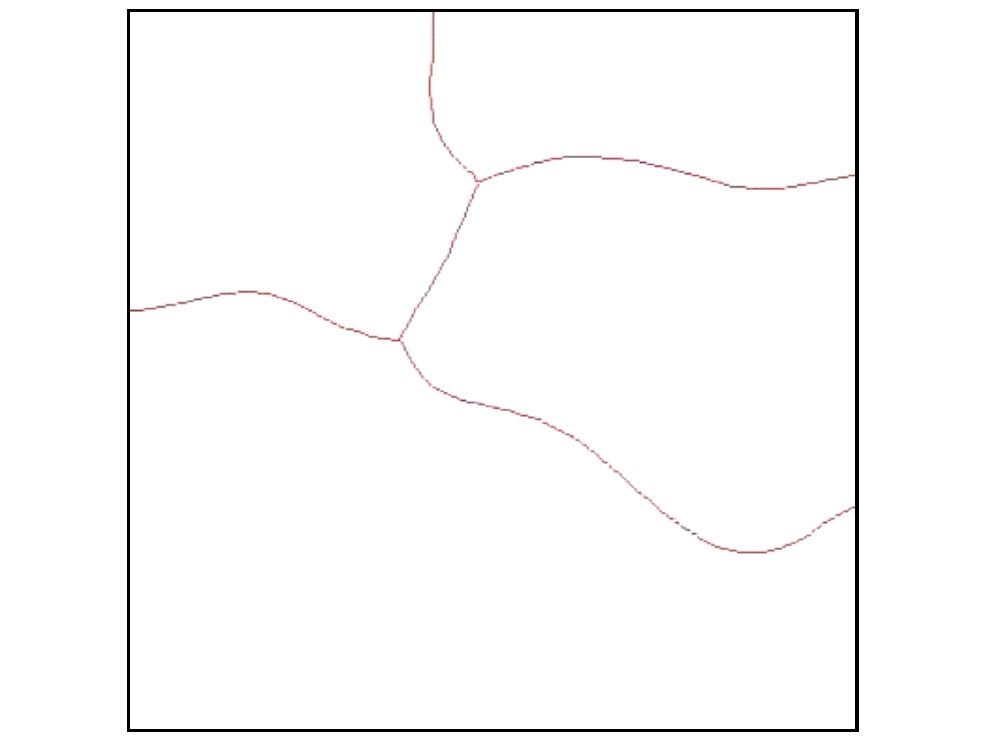} }
		{\includegraphics[scale=0.43,trim=1cm 0cm 1cm 0cm,clip=]{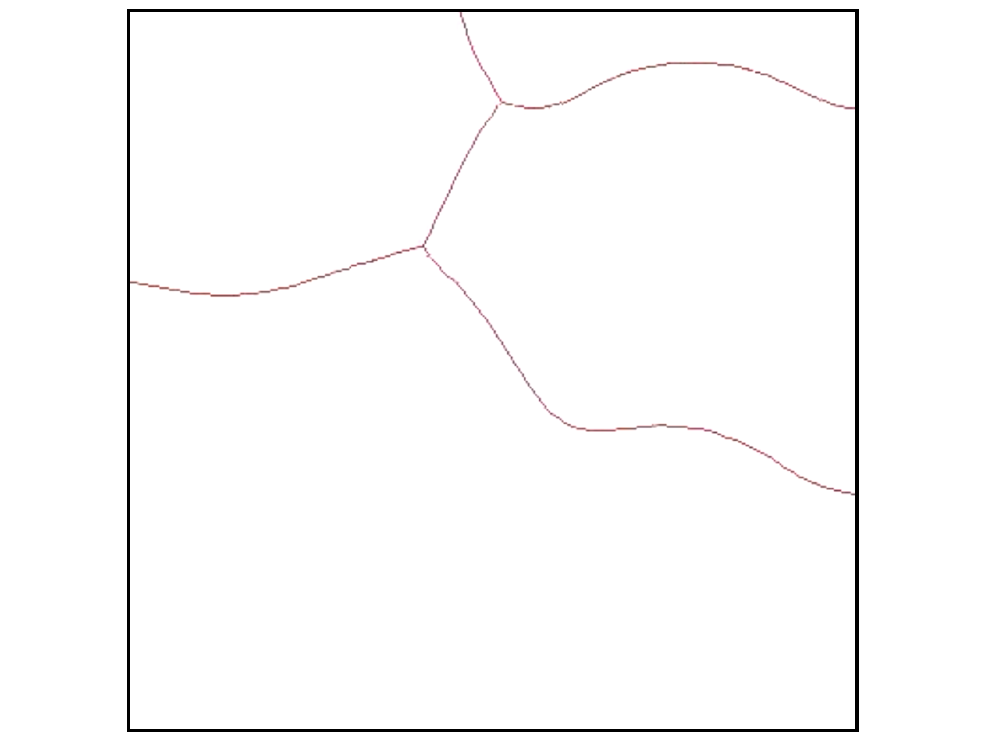} }
	\end{tabular}
   	\caption{Multiphase evolution by HMCF (time is from left to right).}\label{fivebub}
\end{figure}

Our approximation method for computing multiphase interfacial dynamics thus repeats
the PDE step (\ref{vpde}) and the thresholding step (\ref{trunc}).
An example of a numerical result is shown in figure \ref{fivebub}. Here, in contrast to parabolic mean curvature flow, we remark that interfaces are not smoothed, but oscillate. 

\begin{figure}[htbp]
	\begin{center}
		{\includegraphics[scale=0.33,trim=1.5cm 0.0cm 0.5cm 0cm,clip=]{init1.pdf}	}
		{\includegraphics[scale=0.33,trim=1.5cm 0.0cm 0.5cm 0cm,clip=]{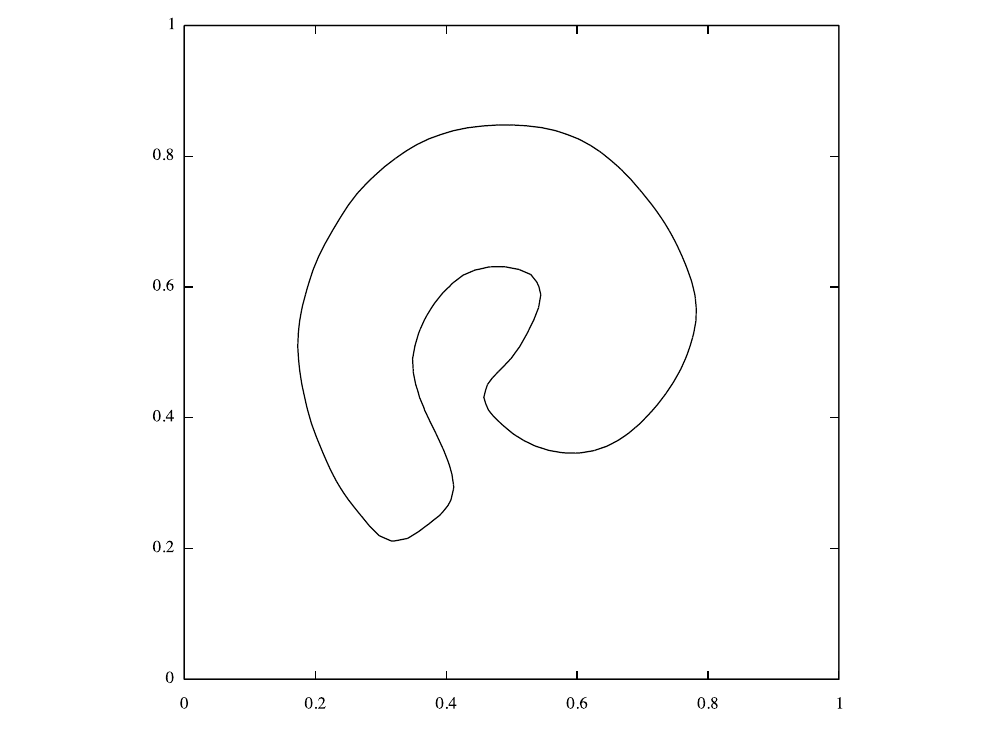}	}
		{\includegraphics[scale=0.33,trim=1.5cm 0.0cm 0.5cm 0cm,clip=]{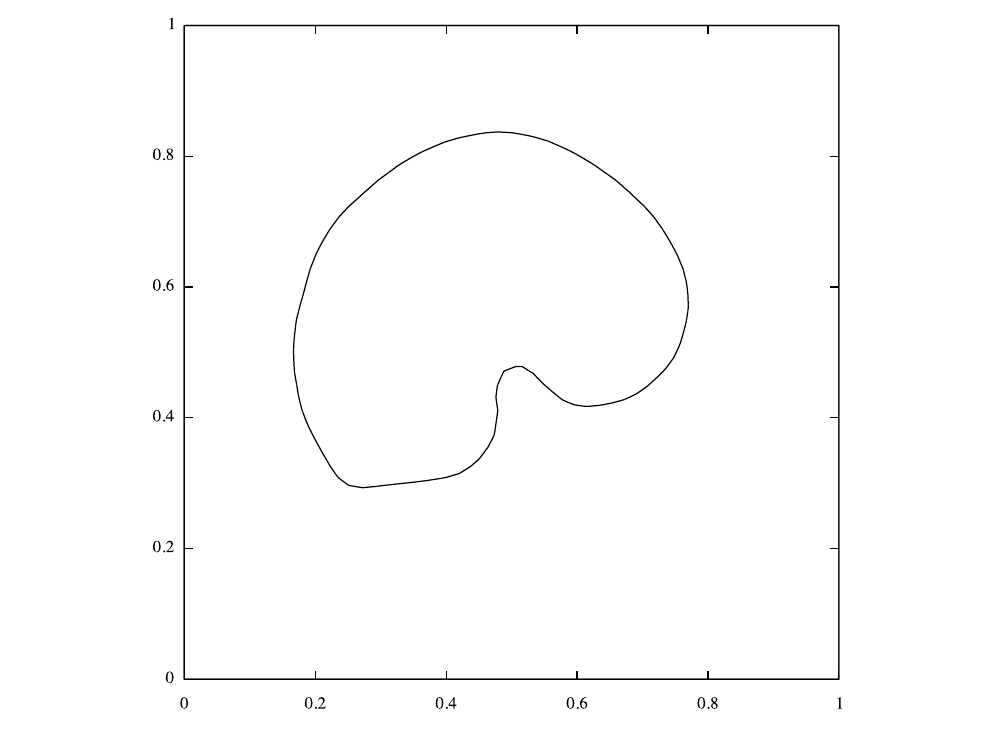}	}
		{\includegraphics[scale=0.33,trim=1.5cm 0.0cm 0.5cm 0cm,clip=]{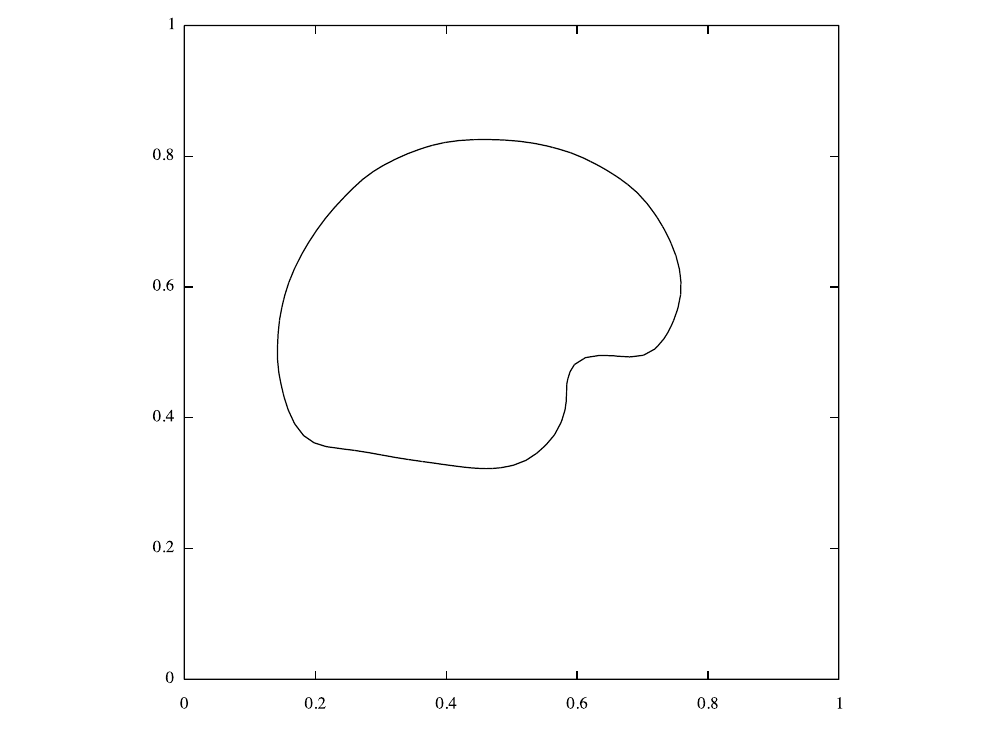}	}
		{\includegraphics[scale=0.33,trim=1.5cm 0.0cm 0.5cm 0cm,clip=]{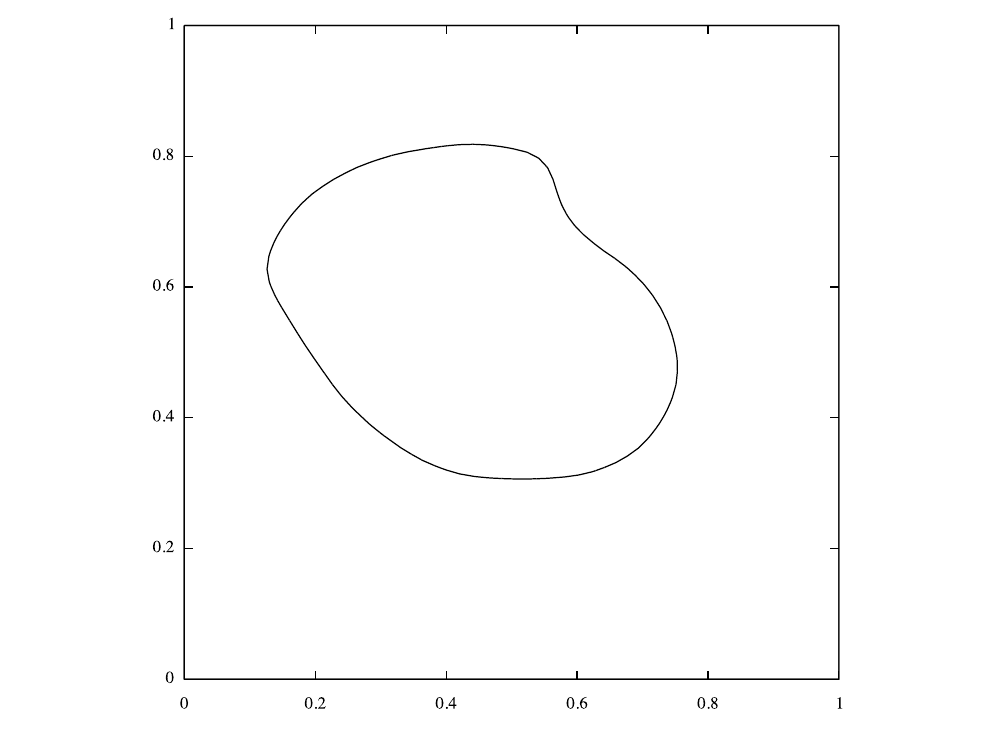}	}
		\caption{Volume preserving motion (time is from left to right).}
		\label{hmcffig2}
	 \end{center}
\end{figure}

As in the parabolic case \cite{[GOS3]}, the method of minimizing movements can be used to inspect volume constrained interfacial motions as well. In particular,  regarding volume-constrained motions, the minimization formulation of our algorithm can allow one to formally include constraints, via penalization. Denoting the prescribed volume of region $P_{i}$ by $A_{i}$, wave-type minimizing movements can be used under the following functional minimization (see, \cite{gs3}):
\begin{equation}\label{vdmf}
	\mathcal{F}_n({\bu})=\IO\left( \frac{|\bu-2\bu_{n-1}+\bu_{n-2}|^2}{2h^2} + \frac{|\nabla \bu|^2}{2}\right){dx}+\frac{1}{\epsilon} \sum_{i=1}^{k} |A_{i}-meas(P_{i}^{\boldsymbol{u}})|^2.
\end{equation}
Here $\bu_{n-1}$ and $\bu_{n-2}$ are created from the wave equation's initial conditions, $\epsilon>0$ is a small penalty parameter and the areas corresponding to ${\boldsymbol{u}}$ are obtained from the sets
\begin{align}
	P_{i}^{\bu} = \{ x \in \Omega ; \;\; {\bu}(x)\cdot {\boldsymbol{p}}_{i} \ge  {\boldsymbol{u}}(x)\cdot {\boldsymbol{p}}_{j} \quad \forall j \}.\notag
\end{align}
Figure \ref{hmcffig2} shows a two phase volume-preserving motion. Contrast to the standard (parabolic) volume-preserving mean curvature flow, whose evolution approaches a circle, the interfacial dynamics considered here oscillate. Using the multiphase formulation of our algorithm, figure \ref{vphmcf} below shows a six phase oscillatory motion. Compuations are performed within the unit square and we remark that volumes are preserved within an absolute error of $10^{-5}$.

\begin{figure}[htbp]
	\begin{center}
		{\includegraphics[scale=0.5,trim=0.0cm 0.0cm 0.5cm 0cm,clip=]{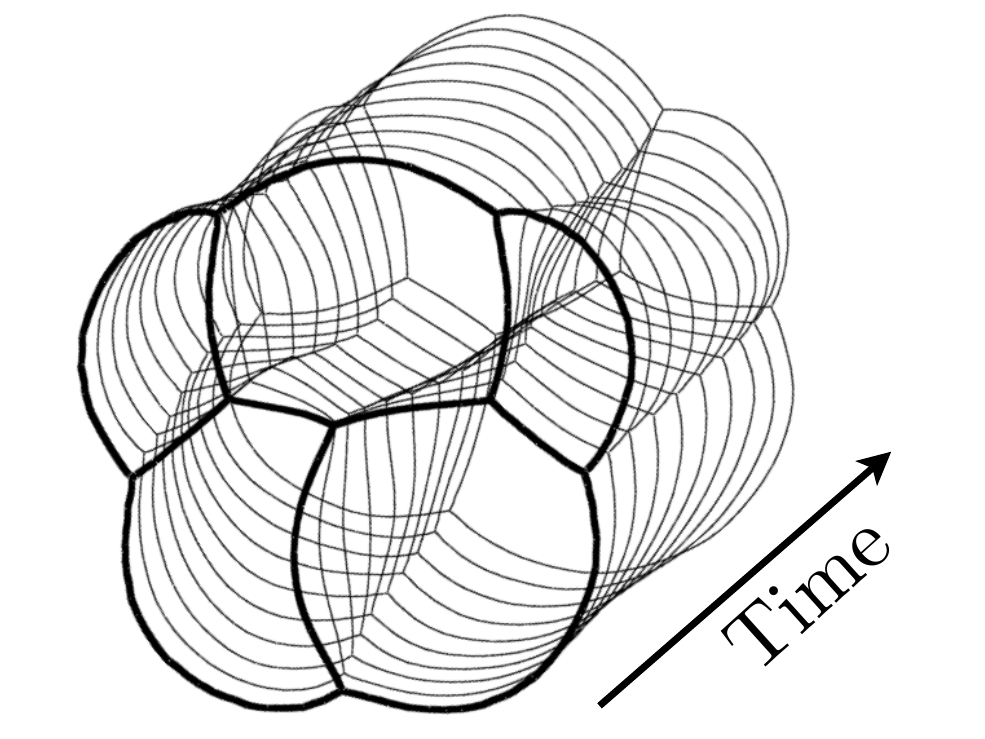}	}
		\caption{Evolution by multiphase volume preserving HMCF.}
		\label{vphmcf}
	 \end{center}
\end{figure}

\section{Conclusion}
We introduced a method for approximating interfacial motion governed by curvature dependent acceleration. Our method is a thresholding algorithm of the BMO-type which, instead of utilizing a diffusion process, thresholds an evolution by the wave equation. We obtained the desired interfacial dynamics by means of a combination of signed distance functions and the convergence order in time of our method was obtained. We also performed a numerical analysis of the algorithm and presented results of its computational application, including investigations of multiphase and volume-preserving motions. The numerical results lead to the expectation of convergence for the volume-preserving threshold dynamical algorithm as well, which we will investigate in the future.

\section{Acknowledgements}
The work of the first author was supported by JSPS Grant Number 25800087 and the second by JSPS Grant Number 26870224.

\end{document}